\documentclass[11pt]{amsart}
\usepackage{amssymb}

\newcommand{\ulambda}{{\boldsymbol{\lambda}}}\newcommand{\umu}{{\boldsymbol{\mu}}}
\newcommand{\unu}{{\boldsymbol{\nu}}}

\newtheorem{thm}{Theorem}[section]

\newtheorem{prop}[thm]{Proposition}

\theoremstyle{definition}
\newtheorem{exmp}[thm]{Example}
\newtheorem{defn}[thm]{Definition}

\theoremstyle{remark}
\newtheorem{rem}[thm]{Remark}

\begin{document}

\title{Crystal graphs of irreducible $\mathcal{U}_v({\widehat{\mathfrak{sl}}_e})$-modules of level two and  Uglov bipartitions}
\author{Nicolas Jacon}

\address{Universit\'e de Franche-Comt\'e,  UFR Sciences et
Techniques, 16 route de Gray, 25 030 Besan\c{c}on, France.}

\email{jacon@math.univ-fcomte.fr}

\date{July, 2006}

\begin{abstract}
We give a simple description of the natural bijection between the set of FLOTW
bipartitions  and the set of Uglov bipartitions (which generalizes the set of
Kleshchev bipartitions). These bipartitions, which  label the crystal graphs of
irreducible $\mathcal{U}_v({\widehat{\mathfrak{sl}}_e})$-modules of level two,
naturally appear in the context of the modular representation theory of Hecke
algebras of type $B_n$.
\end{abstract}

\maketitle

\pagestyle{myheadings}

\markboth{Nicolas Jacon}{Crystal graphs of irreducible
$\mathcal{U}_v({\widehat{\mathfrak{sl}}_e})$-modules of level two}

\maketitle

\section{Introduction}
Let $n>0$ and let $W_n$  be the Weyl group of type $B_n$ with set of simple
reflections $S:=\{t,s_1,\cdots,s_{n-1}\}$ and relations symbolized by the
following braid diagram:
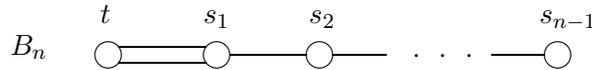
\begin{center}
\begin{picture}(250,30)
\put(  3, 08){$B_n$} \put( 40, 08){\circle{10}} \put( 44, 05){\line(1,0){33}}
\put( 44, 11){\line(1,0){33}} \put( 81, 08){\circle{10}} \put( 86,
08){\line(1,0){29}} \put(120, 08){\circle{10}} \put(125, 08){\line(1,0){20}}
\put(155, 05){$\cdot$} \put(165, 05){$\cdot$} \put(175, 05){$\cdot$} \put(185,
08){\line(1,0){20}} \put(210, 08){\circle{10}} \put( 37, 20){$t$} \put( 76,
20){$s_1$} \put(116, 20){$s_2$} \put(203, 20){$s_{n-1}$}
\end{picture}
\end{center}

Let $k$ be a field and $Q,q\in{k^\times}$. We denote by $H_n:=H_k(W_n,Q,q)$ the
corresponding Iwahori--Hecke algebra. This is an associative unitary
$k$-algebra generated by the elements $T_{s}$ for $s\in{S}$, subject to the
braid relations symbolized by the above diagram and  the relations
$(T_t-Q)(T_t+1)=0$,  $(T_{s_j}-q)(T_{s_j}+1)=0$ for $1\leq j\leq n-1$. When
$H_n$ is semisimple, Tits deformation theorem shows that the simple modules of
this algebra are in natural bijection with the simple modules of the group
algebra $kW_n$. In the non semisimple case, the classification of the simple
$H_n$-modules was achieved by Dipper-James in \cite{DJ},   and Ariki,
Ariki-Mathas in \cite{Ar},\cite{Ar2}, \cite{ArMa} using the theory of canonical
bases and crystal graphs for quantum groups.

  Let  $\mathcal{U}_v({\widehat{\mathfrak{sl}}_e})$ be the quantum group of type
$A^{(1)}_{e-1}$. Then Ariki and Ariki-Mathas have shown that the set of simple
$H_n$--modules $\textrm{Irr}{(H_n)}$ is in natural bijection with the Kashiwara
crystal basis of the irreducible
$\mathcal{U}_v({\widehat{\mathfrak{sl}}_e})$-module with highest weight a sum
of two fundamental weights $\Lambda_{v_0}+\Lambda_{v_1}$ ($0\leq v_0,v_1< e$).
There are several natural ways to obtain a parametrization of this basis,
depending on a choice of integers $s_0$ and $s_1$ in the classes of $v_0$ and
$v_1$ modulo $e$. Hence we obtain several possibilities for labelling the same
set $\textrm{Irr}{(H_n)}$, there are given by a certain class of bipartitions
$\Phi_{e,n}^{(s_0,s_1)}$ named ``Uglov bipartitions". This kind of bipartitions
both generalizes the set of FLOTW bipartitions (which correspond to the case
$0\leq s_0,s_1\leq e$, see \cite{Ja}),  and the set of Kleshchev bipartitions
 (corresponding to the case where $s_0-s_1>n-1-e$, see \cite{Ar2}).

In \cite{GJ}, M. Geck and the author have given an interpretation of this fact
in the context of the representation theory of Hecke algebras. We showed that
each of the parameterizations by $\Phi_{e,n}^{(s_0,s_1)}$ is linked with the
existence of  canonical basic sets and induced the unitriangularity of the
decomposition matrix associated with a  specialization and a choice of a
``weight function" (in the sense of Lusztig \cite[\S 3.1]{Lu}).

In general we only know a recursive definition of the sets of Uglov
bipartitions and a natural problem is to obtain a non recursive (and simple)
characterization of these sets. In the case where $s_0-s_1>n-1-e$ (known as the
``asymptotic case"), this problem has been recently solved by Ariki, Kreiman
and Tsuchioka in \cite{AKT} using results of Littelmann. Our purpose is to
obtain a new characterization of all Uglov bipartitions using the following
facts:
\begin{itemize}
\item in the case where  $0\leq s_0\leq  s_1< e$, we know a simple non recursive
characterization of the set $\Phi_{e,n}^{(s_0,s_1)}$, the FLOTW bipartitions
\cite{FLOTW},
\item if $s_0'\equiv s_0 (\textrm{mod }e)$ and   $s_1'\equiv s_1 (\textrm{mod
}e)$ or if $s_0'\equiv s_1 (\textrm{mod }e)$ and   $s_1'\equiv s_0 (\textrm{mod
}e)$, we have a bijection  between  $\Phi_{e,n}^{(s_0,s_1)}$ and
$\Phi_{e,n}^{(s_0',s_1')}$.
\end{itemize}
Hence, if we know a simple (and non recursive) description of the above
bijection, the desired characterizations of all Uglov bipartitions will follow.
Quite remarkably, the main result of this paper, Theorem \ref{main}, together
with works of Leclerc and Miyachi shows that this bijection is controlled by
the canonical bases of the irreducible
$\mathcal{U}_v({{\mathfrak{sl}}_\infty})$-modules. As a special case, we obtain
a quite simple and new characterization of the set of Kleshchev bipartitions
(but which remains recursive ...) using the notion of symbols.

 The paper will be
organized as follows. The first section  gives a brief exposition of the theory
of crystal graphs and connections with the representation theory of Hecke
algebras. In the second and third section,  our main results are stated and
proved: we study the combinatoric of Uglov bipartitions and we give a
description of the above bijection. This description is largely inspired by the
works of Leclerc and Miyachi. In the last section, we  describe the relations
of our results with these works.

\section{Crystal graphs of $v$-deformed Fock spaces of level $2$}

\subsection{Fock spaces}\label{order} Let $v$ be an indeterminate and let $e$ be a positive integer.
 Let  $\mathfrak{h}$ be a free $\mathbb{Z}$-module with basis
  $\{h_i,\mathfrak{d}\ |\ 0\leq{i}<e\}$ and let $\{\Lambda_i,\delta\ |\ 0\leq{i}<e\}$
  be the dual basis  with respect to the pairing:
$$\langle\ ,\ \rangle :\mathfrak{h}^*\times{ \mathfrak{h}}\to{\mathbb{Z}}$$
such that $\langle\Lambda_i,h_j\rangle=\delta_{ij}$,
$\langle\delta,\mathfrak{d}\rangle=1$ and
$\langle\Lambda_i,\mathfrak{d}\rangle=\langle\delta,h_j\rangle=0$ for
$0\leq{i,j}<e$. The $\Lambda_k$ ($1\leq k\leq e$) are called the {\it
fundamental weights}. The quantum group $\mathcal{U}_v
({\widehat{\mathfrak{sl}}}_e)$ of type $A^{(1)}_{e-1}$  is a unital associative
algebra over $\mathbb{C}(v)$ which is generated by elements $\{e_i,f_i\ |\
i\in{\{0,\cdots,e-1\}}\}$ and $\{k_h\ |\ h\in{\mathfrak{h}}\}$ subject to the
relations described  for example in \cite[chapter 6]{Ma}.

 In
this paper, we want to study the irreducible
$\mathcal{U}_v({\widehat{\mathfrak{sl}}}_e)$-modules with highest weight
$\Lambda$ where $\Lambda$ is a sum of two fundamental weights
$\Lambda_{v_0}+\Lambda_{v_1}$, with $0\leq v_0,v_1<e$. These modules can be
constructed by using the {\it Fock space representation} which we now define.
Let $\Pi_{2,n}$ be the set of bipartitions of rank $n$, that is the set of
$2$-tuples $(\lambda^{(0)},\lambda^{(1)})$ such that $\lambda^{(0)}$ (resp.
$\lambda^{(1)}$) is a partition or rank $a_1$ (resp. $a_2$) with $a_1+a_2=n$.
Let $\boldsymbol{s}=(s_0,s_1)\in{\mathbb{Z}^2}$ be such that $s_0 \equiv v_0\
(\textrm{mod }e)$ and $s_1 \equiv v_1\ (\textrm{mod }e)$ or such that $s_0
\equiv v_1\ (\textrm{mod }e)$ and $s_1 \equiv v_0\ (\textrm{mod }e)$. The {\it
Fock space} (of level $2$) is defined to be the $\mathbb{C}(v)$-vector space
generated by the symbols $|\ulambda,\boldsymbol{s} \rangle$ with $\ulambda\in
\Pi_{2,n}$:
$$\mathfrak{F}^{\boldsymbol{s}}:=\bigoplus_{n\geq 0}\bigoplus_{\ulambda\in{\Pi_{2,n}}}\mathbb{C}(v)
 |\ulambda,\boldsymbol{s} \rangle.$$

 Let us introduce some additional
notations concerning the combinatorics of bipartitions. Let
$\ulambda={(\lambda^{(0)} ,\lambda^{(1)})}$ be a bipartition of rank $n$. The
diagram of $\ulambda$ is the following set:
$$[\ulambda]=\left\{ (a,b,c)\ |\ 0\leq{c}\leq{1},\ 1\leq{b}\leq{\lambda_a^{(c)}}\right\}.$$
The elements of this diagram are called   the {\it nodes} of  $\ulambda$. Let
$\gamma=(a,b,c)$ be a node of  $\ulambda$. The {\it residue} of $\gamma$
associated to $e$ and $({s_0},s_1)$ is the element of $\mathbb{Z}/e\mathbb{Z}$
defined by:
$$\textrm{res}{(\gamma)}\equiv (b-a+s_{c})(\textrm{mod}\ e).$$
If $\gamma $ is a node with residue $i$, we say that  $\gamma$ is an  $i$-node.
Let  $\ulambda$ and  $\umu$ be two bipartitions of rank $n$ and $n+1$ such that
$[\ulambda]\subset{[\umu]}$. There exists a node $\gamma$ such that
$[\umu]=[\ulambda]\cup{\{\gamma\}}$. Then, we denote $[\umu]/[\ulambda]=\gamma$
and  if $\textrm{res}{(\gamma)}=i$, we say that $\gamma$ is an {\it addable}
$i$-node for $\ulambda$ and a {\it removable} $i$-node for  $\umu$. Let
$i\in{\{0,\cdots,e-1\}}$, we introduce a total  order on the set of $i$-nodes
of a bipartition.
 Let $\gamma=(a,b,c)$ and  $\gamma'=(a',b',c')$ be two $i$-nodes of a bipartition. We denote $\gamma<_{(s_0,s_1)}\gamma'$ if:
$$b-a+s_c<b'-a'+s_{c'}\ \textrm{or } \textrm{if}\ b-a+s_c=b'-a'+s_{c'}\textrm{ and }{c'}<{c}.$$
Note that this order strongly depends on the choice of $s_0$ and $s_1$ in the
classes  of $v_0$ and $v_1$ modulo $e$. Note also this the order coincides with
that of \cite{GJ}.

Using this order, it is possible to define an action of
$\mathcal{U}_v({\widehat{\mathfrak{sl}}}_e)$ on the Fock space
$\mathfrak{F}^{\boldsymbol{s}}$   such that $\mathfrak{F}^{\boldsymbol{s}}$
becomes an integrable $\mathcal{U}_v({\widehat{\mathfrak{sl}}}_e)$-module.
Moreover it is known that the submodule $M_{\boldsymbol{s}}$ generated by the
empty bipartition is a highest weight module with weight
$\Lambda_{v_0}+\Lambda_{v_1}$ (\cite{JMMO}  for details). Hence, if ${\bf
s}'=(s_0',s_1')\in{\mathbb{Z}^2}$ is such that $s_0\equiv s_0'(\textrm{mod }e)$
and  $s_1\equiv s_1'(\textrm{mod }e)$ or such that $s_0\equiv s_1'(\textrm{mod
}e)$ and  $s_1\equiv s_0'(\textrm{mod }e)$, then the modules
 $M_{{\bf s}}$ and $M_{{\bf s}'}$ are isomorphic. However,
  it is important to note that  the actions of  $\mathcal{U}_v(\widehat{\mathfrak{sl}}_e)$ on the elements
 of the standard basis $|\ulambda,{\bf s}\rangle$ and
  $|\ulambda,{\bf s}'\rangle$ are different in general.

\begin{rem}
Let $(v_0,v_1)\in{\{0,1,\cdots,e-1\}^2}$. Then it is possible to define another
order on the set of $i$-nodes of a bipartition as follows : we write
$\gamma=(a,b,c)<_{(v_0,v_1)_{+}} \gamma'=(a',b',c')$ if:
$$c'<c \textrm{ or if }c=c' \textrm{ and }a'<a.$$
Note that if we fix a bipartition $\ulambda$ of rank $n$, then the above order
on the $i$-nodes of $\ulambda$ coincides with $<_{(s_0,s_1)}$ in the case where
$s_0\equiv v_0(\textrm{mod }e)$ and  $s_1\equiv v_1(\textrm{mod }e)$ and
$s_0>>s_1$. This order will be referred to the {\it positive asymptotic order}
and this is the one used by Ariki \cite{Ar2} in its determination of the simple
modules for Hecke algebras of type $B_n$.

Similarly, we can define another order on the set of $i$-nodes of a bipartition
as follows: $\gamma=(a,b,c)<_{{(v_0,v_1)}_{-}}\gamma'=(a',b',c')$ if:
$$c'>c \textrm{ or if }c=c' \textrm{ and }a'<a.$$
If we fix a bipartition $\ulambda$ of rank $n$, then  the above order on the
$i$-nodes of $\ulambda$ coincides with $<_{(s_0,s_1)}$  in the case where
$s_0\equiv v_0(\textrm{mod }e)$ and  $s_1\equiv v_1(\textrm{mod }e)$ and
$s_0<<s_1$. This order will be referred to the {\it negative asymptotic order}.

In the two cases, we obtain an action of
$\mathcal{U}_v({\widehat{\mathfrak{sl}}}_e)$ on the space
$\mathfrak{F}^{(v_0,v_1)}$ and the submodules $M^{{+}}_{v_0,v_1}$ and
$M^{{-}}_{v_0,v_1}$ generated by the empty bipartition are both irreducible
highest weight modules with weight $\Lambda_{v_0}+\Lambda_{v_1}$ and they are
isomorphic.

\end{rem}

\subsection{Crystal graph of $M_{{\bf s}}$}\label{crystal} As the modules $M_{{\bf s}}$ are
integrable highest weight modules, the general theory of Kashiwara and Lusztig
provides us with a {\it canonical basis of} $M_{{\bf s}}$. We don't need in
this paper the definition of this basis but by the deep results of Ariki
\cite{Ar}, one of the interest of  this basis is that it provides a way to
compute the decomposition matrices for Hecke algebras of type $B_n$ (see
\cite[Theorem 14.49]{Ar3}). In order to make an efficient use of this, we need
to determine a good parametrization of the canonical basis. This is given by
studying the {\it Kashiwara crystal graph} which we now describe.

Let   $\ulambda$ be a   bipartition and  let $\gamma$ be  an $i$-node of
$\ulambda$,
 we say that  $\gamma$ is a {\it normal} $i$-node of  $\ulambda$ if,
whenever $\eta$ is an $i$-node of $\ulambda$   such that
$\eta>_{(s_0,s_1)}\gamma$, there are more removable $i$-nodes between $\eta$
and $\gamma$ than addable $i$-nodes between $\eta$ and $\gamma$. If $\gamma$ is
the minimal normal $i$-node of $\ulambda$ with respect to $<_{(s_0,s_1)}$, we
say that $\gamma$ is a {\it good} $i$-node.

Following \cite[\S 2]{Ar2}, the normal $i$-nodes of a bipartition $\ulambda$
can be easily obtained using the following process. We first read addable and
removable $i$-nodes of $\ulambda$ in increasing order with respect to
$<_{(s_0,s_1)}$. If we write $A$ for an addable $i$-node and $R$ for a
removable one, we get a sequence of $A$ and $R$. Then we delete $RA$ as many as
possible. The remaining removable $i$-nodes in the sequence are the normal
$i$-nodes and the node corresponding to the leftmost $R$ is a good $i$-node.
\begin{exmp}
Let $e=4$, ${\bf s}=(0,6)$ and $\ulambda=((4,3,1,1),(4))$, The Young diagram of
$\ulambda$ with residues is the following one:
$$  \underline{\lambda}  =\left( \ \begin{tabular}{|c|c|c|c|}
                         \hline
                           0  & 1  & 2  & 3  \\
                          \hline
                            3  & 0 & 1\\
                          \cline{1-3}
                           2  \\
                          \cline{1-1}
                          1 \\
                          \cline{1-1}

                             \end{tabular}\  ,\
                       \begin{tabular}{|c|c|c |c|}
                         \hline
                           2  & 3  & 0 & 1  \\
                          \hline

                             \end{tabular}\
                                                 \right).$$

We have one addable $1$-node $(2,1,1 )$ and three removable $1$-nodes
$(4,1,0)$, $(2,3,0)$ and $(1,4,1)$. We have:
$$(4,1,0)<_{(0,6)} (2,3,0)<_{(0,6)} (2,1,1) <_{(0,6)} (1,4,1),$$
and the associated sequence of removable and addable $1$-nodes is $RRAR$. Hence
$(4,1,0)$ and $(1,4,1)$ are normal $1$-nodes of $\ulambda$ and $(4,1,0)$ is a
good $1$-node for $\ulambda$.
\end{exmp}

 Note that this notion depends on the order $<_{(s_0,s_1)}$ and thus, on the
choice of ${\bf s}$. To define the crystal graph of $M_{{\bf s}}$, we need to
introduce the one of the Fock space $\mathcal{F}^{{\bf s}}$. This graph has
been studied by Jimbo et al. \cite{JMMO}, Foda et al. \cite{FLOTW} and Uglov
\cite{U}. It is given by:
\begin{itemize}
\item vertices: the bipartitions,
\item edges: $\displaystyle{{\ulambda\overset{i}{\rightarrow}{{\umu}}}}$
 if and only if $[\umu]/[\ulambda]$ is a good  $i$-node.
\end{itemize}

Then,  the crystal graph of $M_{{\bf s}}$ is the connected components of that
of $\mathcal{F}^{{\bf s}}$ which contain the vacuum vector
$\boldsymbol{\emptyset}$. The vertices of this graph, which are in natural
bijection with the canonical basis elements of $M_{\bf s}$, are given by the
following class of bipartitions.

\begin{defn}\label{ug} Let ${\bf s}\in{\mathbb{Z}^2}$.
The set of {\it Uglov bipartitions} $\Phi_{e,n}^{\bf s}$ is defined recursively
as follows.
\begin{itemize}
   \item We have $\boldsymbol{\emptyset}:=(\emptyset,\emptyset)\in{\Phi_{e,n}^{\bf s}}$.
    \item If $\ulambda\in\Phi_{e,n}^{\bf s}$ , there exist $i\in{\{0,\cdots,e-1\}}$ and a good $i$-node $\gamma$ such that
    if we remove $\gamma$ from  $\ulambda$, the resulting  bipartition is in $\Phi_{e,n-1}^{\bf s}$.
\end{itemize}
\end{defn}

In the special case where $0\leq s_0 \leq s_1<e$, Foda, Leclerc, Okado, Thibon
and Welsh have given a non recursive parametrization of this kind of
bipartitions.

\begin{prop}[Foda et al. \protect{\cite[Prop. 2.11]{FLOTW}}]\label{FL} Assume that ${\bf
s}:=(s_0,s_1)\in{\mathbb{Z}^2}$ is such that $0\leq s_0 \leq s_1<e$ then
$\ulambda={(\lambda^{(0)},\lambda^{(1)})}$ is in  $\Phi_{e,n}^{\bf s}$  if and
only if:
\begin{enumerate}
\item for all $i=1,2,\cdots$, we have :\begin{align*}
&\lambda_i^{(0)}\geq{\lambda^{(1)}_{i+s_{1}-s_{0}}},\\
&\lambda^{(1)}_i\geq{\lambda^{(0)}_{i+e+s_0-s_{1}}};
\end{align*}
\item  for all  $k>0$, among the residues appearing
at the right ends of the length $k$ rows of   $\ulambda$, at least one element
of  $\{0,1,\cdots,e-1\}$ does not occur.
\end{enumerate}
 Such bipartitions are called FLOTW bipartitions.
\end{prop}

When the condition $0\leq s_0 \leq s_1<e$ isn't satisfied, the above
characterization of Uglov bipartitions is no longer true. Hence, an important
problem would be to obtain a simple description of  $\Phi_{e,n}^{\bf s}$ in all
cases.

Assume that  ${\bf s}:=(s_0,s_1)\in{\mathbb{Z}^2}$ and  ${\bf
s}':=(s'_0,s'_1)\in{\mathbb{Z}^2}$ are such that  $s_0'\equiv s_0 (\textrm{mod
}e)$ and   $s_1'\equiv s_1 (\textrm{mod }e)$ or such that  $s_0'\equiv s_1
(\textrm{mod }e)$ and   $s_1'\equiv s_0 (\textrm{mod }e)$. Then the irreducible
highest weight modules $M_{{\bf s}}$ and  $M_{{\bf s}'}$ are isomorphic and it
implies that the associated Kashiwara crystal graphs are also isomorphic: only
the labelling of the vertices by the sets of Uglov bipartitions changes. Hence,
in these cases, there exists a bijection:
$$\Psi_{(s_0,s_1)}^{(s_0',s_1')}:  \Phi_{e,n}^{(s_0,s_1)} \to \Phi_{e,n}^{(s_0',s_1')}.$$
 This bijection may be obtained by following a sequence of arrows back to the
 empty bipartition in the crystal graph of  $M_{{\bf s}}$ and then applying
 the reversed sequence to the empty bipartition of  $M_{{\bf s}'}$. In other
 words, the bijection is obtained recursively as follows. We put
 $\Psi_{(s_0,s_1)}^{(s_0',s_1')}({\boldsymbol{\emptyset}})={\boldsymbol{ \emptyset}}$. Assume
 that we know $\Psi_{(s_0,s_1)}^{(s_0',s_1')}:  \Phi_{e,n-1}^{(s_0,s_1)} \to \Phi_{e,n-1}^{(s_0',s_1')}.$
 Let $\ulambda\in{\Phi_{e,n}^{(s_0,s_1)}
 }$. Then, there exist $i\in{\{0,\cdots,e-1\}}$ and a good $i$-node $\gamma$ with respect to $<_{(s_0,s_1)}$ such that
    if we remove $\gamma$ from  $\ulambda$, the resulting  bipartition $\ulambda'$ is in $\Phi_{e,n-1}^{\bf
    s}$. Let $\umu':=\Psi_{(s_0,s_1)}^{(s_0',s_1')}(\ulambda')$. Then there
    exist an $i$-node $\gamma'$ and a bipartition $\umu$ such that  $[\umu]=[\umu']\cup{\{\gamma'
    \}}$ and such that $\gamma'$ is a good $i$-node for $\umu$ with respect to
    $<_{(s_0',s_1')}$. Then, we  put
    $\Psi_{(s_0,s_1)}^{(s_0',s_1')}(\ulambda)=\umu$.

\begin{rem}\label{asyp}
Let $(v_0,v_1)\in{\{0,1,\cdots,e-1\}}^2$. Then the crystal associated to the
modules $M^{{+}}_{v_0,v_1}$ and $M^{{-}}_{v_0,v_1}$ may be obtained by the same
way as Def. \ref{ug} using the order $<_{(v_0,v_1)_{+}}$ and
$<_{(v_0,v_1)_{-}}$. The bipartitions which label the vertices of the crystal
graph are respectively called the {\it positive  Kleshchev bipartitions} and
the {\it negative  Kleshchev bipartitions}. They are denoted by
$\Phi_{e,n}^{(v_0,v_1)_{+}}$ and $\Phi_{e,n}^{(v_0,v_1)_{-}}$.

 Let ${\bf s}:=(s_0,s_1)\in{\mathbb{Z}^2}$ be such that  $s_0\equiv v_0 (\textrm{mod
}e)$ and   $s_1\equiv v_1 (\textrm{mod }e)$. Then the irreducible highest
weight modules $M_{{\bf s}}$, $M^{{+}}_{v_0,v_1}$ and $M^{{-}}_{v_0,v_1}$ are
isomorphic and we also obtain bijections:
$$\Psi_{(s_0,s_1)}^{(v_0,v_1)_{-}   }:  \Phi_{e,n}^{(s_0,s_1)} \to \Phi_{e,n}^{(v_0,v_1)_{-}},$$
$$\Psi_{(s_0,s_1)}^{(v_0,v_1)_{+}}:  \Phi_{e,n}^{(s_0,s_1)} \to \Phi_{e,n}^{(v_0,v_1)_{+}},$$
$$\Psi_{ (v_0,v_1)_{+}   }^{(v_0,v_1)_{-}}:  \Phi_{e,n}^{(v_0,v_1)_{+} } \to \Phi_{e,n}^{(v_0,v_1)_{-}}.$$
Note that we also have bijections $\Psi_{(s_0,s_1)}^{(v_1,v_0)_{-} }$,
$\Psi_{(s_0,s_1)}^{(v_1,v_0)_{+}}$ and $\Psi_{ (v_0,v_1)_{+}
}^{(v_1,v_0)_{-}}$.  By the definitions of the order $<_{(v_0,v_1)_{+}}$ and
$<_{(v_1,v_0)_{-}}$ and the definition of good nodes, it is clear that the last
bijection is given by
$\Psi_{(v_0,v_1)_{+}}^{(v_1,v_0)_{-}}(\lambda^{(0)},\lambda^{(1)})=(\lambda^{(1)},\lambda^{(0)})$
for all $(\lambda^{(0)},\lambda^{(1)})\in{ \Phi_{e,n}^{(v_0,v_1)_{+}}}$.

\end{rem}

Now it is natural to try to obtain a more efficient description of these
bijections. This is also motivated by the following results.

\subsection{Hecke algebras of type $B_n$} One of the motivations for studying
the class of Uglov bipartitions is provided by the study of the modular
representations of Hecke algebras of type $B_n$. We briefly sketch this
application in this subsection.

Let $W_n$ be the Weyl group of type $B_n$, let $(a,b)\in{\mathbb{N}_{>0}^2}$
and  $\zeta_l:=\textrm{exp}(\frac{2i\pi}{l})$. Let
$H_n:=H_k(W_n,\zeta_l^b,\zeta_l^a)$ be the Hecke algebra with parameters
$Q:=\zeta_l^b$ and $q:=\zeta_l^a$  defined over the field of complex numbers as
it is defined in the introduction. In this case, the algebra $H_n$ is non
semisimple in general and one of the main problem is to determine a
parametrization of its simple modules and to compute the associated
decomposition matrix. An approach to solve this problem has been given by Geck
\cite{Gm} and Geck-Rouquier \cite{GR}. This approach which is closely related
to the existence of Kazhdan-Lusztig theory shows the existence of  ``canonical
sets" of bipartitions which are in natural bijection with the set
$\textrm{Irr}(H_n)$. These sets are called ``canonical basic sets" and they
also show the unitriangularity of the decomposition matrix of $H_n$ (for a good
order on the rows provided by Lusztig $a$-function). A complete survey of this
theory can be found in \cite{Gm} (see also \cite{Gp} for further applications).
Now, \cite[Theorem 5.4]{GJ} shows  that these canonical basic sets are
precisely given by the Uglov bipartitions.
\begin{thm}[Geck-Jacon \cite{GJ}] Let $H_n:=H_k(W_n,\zeta_l^b,\zeta_l^a)$ be the Hecke algebra with parameters
$Q:=\zeta_l^b$ and $q:=\zeta_l^a$ where $(a,b)\in{\mathbb{N}_{>0}^2}$. Let
$d\in{\mathbb{Z}}$ be such that $\zeta_l^b=-\zeta_l^{a.d}$. Let $e\geq 2$ be
the multiplicative order of $q$ and let $p\in\mathbb{Z}$ be such that:
$$d+p e <\frac{b}{a}<d+(p+1)e.$$
Then the set $\mathcal{B}= \Phi_{e,n}^{(d+pe,0)}$ is a canonical basic set in
the sense of \cite[Def. 2.4]{GJ} and it is in natural bijection with
$\operatorname{Irr}(H_n)$.

\end{thm}

Thus it could be interesting to obtain another characterization of the set of
Uglov bipartitions.

\section{First results}
In this section, we show that the characterization of the map
$\Psi_{(s_0,s_1)}^{(s_0,s_1+e)}$ in the case where $s_0\leq s_1$ is sufficient
to obtain a characterization of the maps $\Psi_{(s_0,s_1)}^{(s_0',s_1')}$ in
all cases.

\subsection{Particular cases    } the following proposition gives the explicit
description of the map $\Psi_{(s_0 ',s_1 ')}^{(s_0,s_1)}$ in particular cases.
\begin{prop}\label{pr1}
Let $(s_0,s_1)\in{\mathbb{N}^2}$ and let $e$ be a positive integer such that
$e>1$.
\begin{enumerate}
\item Let $t\in{\mathbb{Z}}$, then for all $\ulambda\in{\Phi_{e,n}^{(s_0,s_1)}}$
we have $\Psi_{(s_0,s_1)}^{(s_0+te,s_1+te)} (\ulambda)=\ulambda$. Hence we have
$$\Phi_{e,n}^{(s_0,s_1)} = \Phi_{e,n}^{(s_0+te,s_1+te)}.$$
\item For all
$\ulambda=(\lambda^{(0)},\lambda^{(1)})\in{\Phi_{e,n}^{(s_0,s_1)}}$ we have
$\Psi_{(s_0,s_1)}^{(s_1,s_0+e)}  (\lambda^{(0)},\lambda^{(1)}  )=
(\lambda^{(1)},\lambda^{(0)}) $. Hence we have
$$\Phi_{e,n}^{(s_1,s_0+e)} = \left\{   \ulambda=(\lambda^{(0)},\lambda^{(1)})\in{\Pi_{2,n}}\ |\  (\lambda^{(1)},\lambda^{(0)})\in{
 \Phi_{e,n}^{(s_0,s_1)}}  \right\}.$$

\end{enumerate}

\end{prop}

\begin{proof}
The first assertion is clear as the order associated to $(s_0,s_1)$ and
$(s_0+te,s_1+te)$  on the set of $i$-nodes of a bipartition is the same in both
cases.

We prove $(2)$ by induction on the rank $n$. If $n=0$, then the result is
clear. Assume that $n>0$. Let
$\ulambda=(\lambda^{(0)},\lambda^{(1)})\in{\Phi_{e,n}^{(s_0,s_1)}}$ and let
$\gamma=(a,b,c)$ be a good $i$-node of $\ulambda$. We must show that
$\gamma'=(a,b,c+1(\textrm{mod }2  ))$ is a good $i$-node for
$(\lambda^{(1)},\lambda^{(0)})$  for the order induced by $(s_1,s_0+e)$ and the
result will follow by induction. To do this, by the definition of good nodes \S
in \ref{crystal}, it is enough to show the following property: let
$i\in{\{0,1,...,e-1\}}$ then $\gamma_1=(a_1,b_1,c_1)$ is an $i$-node in $
(\lambda^{(0)},\lambda^{(1)})$ such that $\gamma>_{(s_0,s_1)}\gamma_1$ if and
only if $\gamma_1 '=(a_1,b_1,c_1 +1(\textrm{mod }2  ))$ is an $i$-node in $
(\lambda^{(1)},\lambda^{(0)})$ such that $\gamma'>_{(s_1,s_0+e)}\gamma_1 '$.
 We first assume  that  $\gamma>_{(s_0,s_1)}\gamma_1$ and we show $\gamma'>_{(s_1,s_0+e)}\gamma_1
 '$. Note that as $\gamma$ and $\gamma_1$ have the same residue modulo $e$,
 there exists $t\in{\mathbb{Z}}$ such that $b-a+s_c=b_1-a_1+s_{c_1}+t e $.
\begin{itemize}
\item If $c=c_1$ then it is clear that $\gamma'>_{(s_1,s_0+e)}\gamma_1$.
\item if $c=0$ and $c_1=1$ then we have $t\geq 0$. Hence
 $b-a+s_0\geq b_1-a_1+s_{1}$ thus $b-a+s_0+e> b_1-a_1+s_{1}$  and   $\gamma'>_{(s_1,s_0 +e)}\gamma_1
 '$.

\item if $c=1$ and $c_1=0$ then we have $t> 0$. Hence we have
 $b-a+s_1\geq  b_1-a_1+s_{0}+e$. If $t>1$ then we have $b-a+s_1> b_1-a_1+s_{0}+e
 $ and $\gamma'>_{(s_1,s_0+e)}\gamma_1
 '$. If $t=1$ then  we have $b-a+s_1= b_1-a_1+s_{0}+e $ and $\gamma'>_{(s_1,s_0+e)}\gamma_1
 '$ because $\gamma'$ is in the first component of $(\lambda^{(1)},\lambda^{(0)})$.
\end{itemize}
Assume now that  $\gamma'>_{(s_1+e,s_0)}\gamma_1
 '$ then by the above argument  $\gamma>_{(s_0+e,s_1+e)}\gamma_1$ and we
 conclude using $(1)$.
\end{proof}

The following proposition deals with the characterization of the maps
$\Psi_{(s_0,s_1)}^{(v_0,v_1)_{-}}$ and $\Psi_{(s_0,s_1)}^{(v_0,v_1)_{+}}$.

\begin{prop}\label{pr2}
Let $(s_0,s_1)\in{\mathbb{Z}^2}$ and let $(v_0,v_1)\in{\{0,1,...,e-1\}}^2$ be
such that $v_0\equiv s_0 (\textrm{mod }e)$ and $v_1\equiv s_1 (\textrm{mod
}e)$.
\begin{enumerate}
\item  if $s_1-s_0>n-1$ then for all
$\ulambda\in{\Phi_{e,n}^{(s_0,s_1)}}$ we have $\Psi_{(s_0,s_1)}^{(v_0,v_1)_{-}}
(\ulambda)=\ulambda$. Hence we have
$$\Phi_{e,n}^{(s_0,s_1)} = \Phi_{e,n}^{(v_0,v_1)_{-} }.$$
\item If $s_0-s_1>n-1-e$ then for all $\ulambda\in{\Phi_{e,n}^{(s_0,s_1)}}$ we have
$\Psi_{(s_0,s_1)}^{(v_0,v_1)_{+}} (\ulambda)=\ulambda$. Hence we have
$$\Phi_{e,n}^{(s_0,s_1)} = \Phi_{e,n}^{(v_0,v_1)_{+} }.$$
\end{enumerate}

\end{prop}

\begin{proof}
We prove $(1)$. Let $(s_0,s_1)\in{\mathbb{Z}^2}$ be such that $s_1-s_0>n-1$ and
let $\ulambda\in{\Phi_{e,n}^{(s_0,s_1)}}$. Let $\gamma=(a,b,c)$ be a removable
$i$-node of $\ulambda$ and let $\gamma'=(a',b',c')$ be an addable or removable
$i$-node of $\ulambda$. We show that $\gamma<_{(s_0,s_1)}\gamma'$ if and only
if $\gamma<_{(v_0,v_1)_{-}}\gamma'$ and the result will follow by induction and
by the definition of good $i$-node as in the proof of the previous Proposition.
Assume first that $\gamma<_{(s_0,s_1)}\gamma'$.  If $c=c'$ then the result is
clear. So assume that $c\neq c'$. If $c'=1$ and $c=0$ then we have
$\gamma<_{(v_0,v_1)_{-}}\gamma'$. Assume that $c=1$ and $c'=0$. As $\gamma$ and
$\gamma'$ have the same residue modulo $e$, there exists $t\in{\mathbb{Z}}$
such that $b-a+s_1=b'-a'+s_{0}+t e $. As $\gamma<_{(s_0,s_1)}\gamma'$, we have
$t\leq 0$. Hence we have :
 $$b-a-(b'-a')\leq (s_0-s_{1})<1-n$$
This is impossible. Indeed, as $\ulambda$ is a bipartition of rank $n$, we must
have:
$$|b'-a'-(b-a)|\leq n-1.$$
Assume now that $\gamma<_{(v_0,v_1)_{-}}\gamma'$. If  $c=c'$ then
$\gamma<_{(s_0,s_1)}\gamma'$.  If otherwise, we have $c'=1$ and $c=0$ then
$b'-a'+s_1-(b-a+s_0)\geq 1-n+s_1-s_0>0$ and we conclude that
$\gamma<_{(s_0,s_1)}\gamma'$. Hence the first assertion is proved. $(2)$
follows by using Prop. \ref{pr1} $(2)$ and Remark \ref{asyp}.

\end{proof}

\subsection{The map $\Psi_{(s_0,s_1)}^{(s_0,s_1+e)}$}
In this subsection, we show that it is enough to characterize
$\Psi_{(s_0,s_1)}^{(s_0,s_1+e)}$ in the case where $s_0\leq s_1$ to
characterize $\Psi_{(s_0,s_1)}^{(s_0',s_1')}$ in all cases.

So, let's assume that we know $\Psi_{(s_0,s_1)}^{(s_0,s_1+e)}$ and its reversed
map if $s_0\leq s_1$.

Let $(u_0,u_1)\in{\mathbb{Z}^2}$. By Prop \ref{pr1} $(2)$, we can assume that
$0\leq u_0\leq u_1<e$. Then,  we have a characterization of all the following
maps:

$$\begin{array}{ccccccccc}

&\Phi_{e,n}^{(u_0,u_1)} &\xrightarrow[\Psi_{(u_0,u_1)}^{(u_0,u_1+e)}]{}
&\Phi_{e,n}^{(u_0,u_1+e)}&
 \xrightarrow[\Psi_{(u_0,u_1+e)}^{(u_0,u_1+2e)}]{} &
\cdots & & \\
&  \cdots  & \xrightarrow[\Psi_{(u_0,u_1+(t-1)e)}^{(u_0,u_1+te)}]{}   &
\Phi_{e,n}^{(u_0,u_1+te)}&\xrightarrow[\Psi_{(u_0,u_1+te)}^{(u_0,u_1+(t+1)e )
}]{} & \Phi_{e,n}^{(u_0,u_1)_{-} }.\end{array}$$
 where $t$ is such that $(t-1)
e>n-1$.
  Now, by Prop \ref{pr1} $(1)$, for all $s\in{\mathbb{N}}$ we have
$\Phi_{e,n}^{(u_0,u_1+se)}=\Phi_{e,n}^{(u_0-se,u_1)}$ and
$\Psi_{(u_0,u_1+se)}^{(u_0-se,u_1)}$ is the identity. Hence all the following
maps are known:
$$\begin{array}{ccccccccc}
& \Phi_{e,n}^{(u_0,u_1)} & \xrightarrow[\Psi_{(u_0,u_1)}^{(u_0-e,u_1)}]{}
&\Phi_{e,n}^{(u_0-e,u_1)} & \xrightarrow[\Psi_{(u_0-e,u_1)}^{(u_0-2e,u_1)}]{} &
\cdots &&\\
& \cdots & \xrightarrow[\Psi_{(u_0-(t-1)e,u_1)}^{(u_0-te,u_1)}]{}
&\Phi_{e,n}^{(u_0-t e,u_1)}&\xrightarrow[\Psi_{(u_0-t e,u_1)}^{(u_0-(t+1)
e,u_1)
 }]{} & \Phi_{e,n}^{(u_0,u_1)_{-} }.\end{array}$$

As we have $0\leq u_0 \leq u_1<e$, we have $u_1\leq u_0+e$. Hence, we have a
characterization of the following maps:
$$\begin{array}{ccccccccc}
& \Phi_{e,n}^{(u_1,u_0+e)}&  \xrightarrow[\Psi_{(u_1,u_0 +e)}^{(u_1,u_1+2e)}]{}
& \Phi_{e,n}^{(u_1,u_0+2e)} &
 \xrightarrow[\Psi_{(u_1,u_0+2e)}^{(u_1,u_0+3e)}]{} &
\cdots && \\
& \cdots & \xrightarrow[\Psi_{(u_1,u_0+te)}^{(u_1,u_0+(t+1)e)
 }]{}
&\Phi_{e,n}^{(u_1,u_0+(t+1)e)}&
\xrightarrow[\Psi_{(u_1,u_0+(t+1)e)}^{(u_1,u_0+(t+2)e)
 }]{} & \Phi_{e,n}^{(u_1,u_0)_{-} }.\end{array}$$
Hence by Prop \ref{pr1} $(2)$, we have a characterization of the following
maps:
$$\begin{array}{ccccccccc}
 & \Phi_{e,n}^{(u_0,u_1)}  & \xrightarrow[\Psi_{(u_0,u_1)}^{(u_0+e,u_1)}]{}  &  \Phi_{e,n}^{(u_0+e,u_1)}
 & \xrightarrow[\Psi_{(u_0+e,u_1)}^{(u_0+2e,u_1)}]{} &
\cdots &&\\
& \cdots
  & \xrightarrow[\Psi_{(u_0+(t-1)e,u_1)}^{(u_0+te,u_1
 )}]{}  &\Phi_{e,n}^{(u_0+te,u_1)} & \xrightarrow[\Psi_{(u_0+te,u_1)}^{(u_0+(t+1)e,u_1
 )}]{} &  \Phi_{e,n}^{(u_0,u_1)_{+} }.\end{array}$$
By Prop \ref{pr1} $(1)$, for all $s\in{\mathbb{N}}$ we have
$\Phi_{e,n}^{(u_0+se,u_1)}=\Phi_{e,n}^{(u_0,u_1-se)}$ and
$\Psi_{(u_0,u_1+se)}^{(u_0-se,u_1)}$ is the identity. Hence all the following
maps are known:
$$\begin{array}{ccccccccc}
& \Phi_{e,n}^{(u_0,u_1)}  & \xrightarrow[\Psi_{(u_0,u_1)}^{(u_0,u_1-e)}]{}&
\Phi_{e,n}^{(u_0,u_1-e)}&
 \xrightarrow[\Psi_{(u_0,u_1-e)}^{(u_0,u_1-2e)}]{} &
\cdots & & \\
& \cdots & \xrightarrow[\Psi_{(u_0,u_1-(t-1)e)}^{(u_0,u_1-t+e
 )}]{}
&\Phi_{e,n}^{(u_0,u_1-te)}& \xrightarrow[\Psi_{(u_0,u_1-te)}^{(u_0,u_1-(t+1)e
 )}]{} & \Phi_{e,n}^{(u_0,u_1)_{+} }.\end{array}$$
Thus, we conclude that the characterization of $\Psi_{(s_0,s_1)}^{(s_0,s_1+e)}$
in the case where $s_0\leq s_1$ yields a characterization of
$\Psi_{(s_0,s_1)}^{(s_0',s_1')}$ in all cases.

\section{Characterization of the map $\Psi_{(s_0,s_1)}^{(s_0,s_1+e)}$}

\subsection{Properties of Uglov bipartitions} We begin with a general result on the set of
Uglov bipartitions. This will be useful for the proof of the main result.
\begin{prop}\label{std} Let ${\bf s}:=(s_0,s_1)\in{\mathbb{Z}^2}$ and
assume that $s_1\geq s_0$. Let $\ulambda\in{\Phi_{e,n}^{\bf s}}$ then
$\ulambda\in{\Phi_{f,n}^{\bf s}}$ where $f>\operatorname{Max}(s_0+n,s_1+n)$.
Hence, for all $i=1,2,\cdots$ we have :
$$\lambda_i^{(0)}\geq{\lambda^{(1)}_{i+s_{1}-s_{0}}}.$$
\end{prop}
\begin{proof}
This is proved by induction on $n$. If $n=0$, the result is trivial. Let $n>0$
and let  $\ulambda\in{\Phi_{e,n}^{\bf s}}$. Then by the definition of Uglov
bipartitions, there exists a good $i$-node $\eta=(a,b,c)$ such that if we
remove $\eta$ from  $\ulambda$, the resulting  bipartition is in
$\Phi_{e,n-1}^{\bf s}$. We have $\lambda^{(c)}_a-a+s_c\equiv i (\textrm{mod
}e)$. Now, we have two cases to consider:
\begin{itemize}
\item If there is no addable node $\eta'=(a',b',c')$ such that $\lambda^{(c')}_{a'}-a'+s_{c'}
= \lambda^{(c)}_a-a+s_c$, as   $f>\operatorname{Max}(s_0+n,s_1+n)$,  there is
no addable node such that $\lambda^{(c')}_{a'}-a'+s_{c'} \equiv
\lambda^{(c)}_a-a+s_c (\textrm{mod }f)$. It implies that $\eta$ is a normal
node for the order induced by ${\bf s}$ and $f$. If there is no removable node
$\eta'=(a',b',c')$ such that $\lambda^{(c')}_{a'}-a'+s_{c'} =
\lambda^{(c)}_a-a+s_c$ then this is a good node for the order induced by ${\bf
s}$ and $f$. If otherwise, as $\eta$ is a good node for the order induced by
${\bf s}$ and $e$, we must have $c'<c$. We conclude that $\eta$ is a good
$i$-node for the order induced by ${\bf s}$ and $f$.
\item  If there is an addable node $\eta'=(a',b',c')$ such that $\lambda^{(c')}_{a'}-a'+s_{c'}
= \lambda^{(c)}_a-a+s_c$, then, as $\eta$ is a good $i$-node for the order
induced by ${\bf s}$ and $e$, we must have $c'>c$ (if otherwise, we have
$\eta'>_{(s_0,s_1)}\eta$ and no removable $i$-node between these two
$i$-nodes). $\eta'$ is the only addable node which has the same residue as
$\eta'$ modulo $f$. Moreover, in this case, there is no removable node
$\eta'=(a',b',c')$ such that $\lambda^{(c')}_{a'}-a'+s_{c'} =
\lambda^{(c)}_a-a+s_c$ and thus, such that $\lambda^{(c')}_{a'}-a'+s_{c'}
\equiv \lambda^{(c)}_a-a+s_c (\textrm{mod }f)$. Hence $\eta$ must be a good
$i$-node for the order induced by ${\bf s}$ and $f$.
\end{itemize}
Thus, the first part of the proposition  follows by induction.

Now, as $f>\operatorname{Max}(s_0+n,s_1+n)$, the elements of $\Phi_{f,n}^{\bf
s}$ are FLOTW bipartitions. Hence, we can use the characterization of Prop.
\ref{FL} to get the second part of the proposition.
\end{proof}

\subsection{Symbol of a bipartition}\label{symbol} Let
$\boldsymbol{s}:=(s_0,s_1)\in{\mathbb{N}^2}$ be such that $s_0\leq s_1$
 and let  ${\ulambda}:=(\lambda^{(0)},\lambda^{(1)})$ be a bipartition of
 rank $n\geq 0$.
 Assume that
 $\lambda^{(0)}=(\lambda^{(0)}_1,\lambda^{(0)}_2,\cdots,\lambda^{(0)}_{r_0})$ and
 $\lambda^{(1)}=(\lambda^{(1)}_1,\lambda^{(1)}_2,\cdots,\lambda^{(1)}_{r_1})$
 (where $\lambda^{(0)}_1 \geq \lambda^{(0)}_2\geq \cdots\geq \lambda^{(0)}_{r_0}$
 and $\lambda^{(1)}_1 \geq \lambda^{(1)}_2\geq \cdots\geq \lambda^{(1)}_{r_1}$).
Let $m\in{\mathbb{N}}$ be such that $m> \textrm{Max}(r_0-s_0,r_1-s_1)$. We
define the following numbers which depends on ${\ulambda}$, ${\bf s}$  and $m$:
\begin{itemize}
\item for $i=1,\cdots,m+s_0$, we put $\beta^{(0)}_i=\lambda_i^{(0)}-i+s_0+m$,
\item for $j=1,\cdots,m+s_1$, we put $\beta^{(1)}_j=\lambda^{(1)}_j-j+s_1+m$.
\end{itemize}
where we put $\lambda^{(0)}_k:=0$ (resp.  $\lambda^{(1)}_k:=0$) if $k>r_0$
(resp. $k>r_1$). We have
$\beta^{(1)}_1>\beta^{(1)}_2>\cdots>\beta^{(1)}_{m+s_1}\geq 0$ and
$\beta^{(0)}_1>\beta^{(0)}_2>\cdots>\beta^{(0)}_{m+s_0}\geq 0$. Then, the ${\bf
s}$-symbol $S_{\bf s}({\ulambda})$ of ${\ulambda}$ is define to be the pair of
these two partitions. This is written as follows:
 $$\left(\begin{array}{ccccc}
 \beta^{(1)}_{m+s_1}  & \beta^{(1)}_{m+s_1-1} & \cdots &  \cdots & \beta^{(1)}_1 \\
      \beta^{(0)}_{m+s_0} & \beta^{(0)}_{m+s_0-1} & ...&  \beta^{(0)}_1 &
\end{array}
\right)$$ On the other hand, given a ${\bf s}$-symbol $S_{\bf s}$, it is easy
to get the bipartition ${\ulambda}$ such that $S_{\bf s}=S_{\bf s}(\ulambda)$.

By Proposition \ref{std}, note that if $\ulambda$ is in  $\Phi_{e,n}^{\bf s}$,
the ${\bf s}$-symbol $S_{\bf s}({\ulambda})$ has the property that
$\beta^{(1)}_i\leq \beta^{(0)}_i$ for $i=1,\cdots,m+s_0$. Such symbols are
called {\it standard} in \cite{LM}.

We will now define a map from the set of Uglov bipartitions $\Phi_{e,n}^{\bf
s}$ to the set of bipartitions of rank $n$ using this notion of ${\bf
s}$-symbol. Let ${\ulambda}:=(\lambda^{(0)},\lambda^{(1)})\in\Phi_{e,n}^{\bf
s}$  and let $S_{\bf s}({\ulambda})=\left(\begin{array}{cc}
\beta^{(1)}\\
\beta^{(0)}
\end{array}
\right) $ be its ${\boldsymbol
 s}$-symbol. Following \cite[\S 2.5]{LM}, we first define an injective map $\theta:\beta^{(0)}\to
 \beta^{(1)}$ such that $\theta(\beta^{(0)}_j)\leq \beta^{(0)}_j$ for all $j\in{\{1,\cdots,m+s_0\}}$ as follows.
 \begin{itemize}
 \item Let $\beta^{(1)}_i$ be the maximal element of $\beta^{(1)}$ such that
 $\beta^{(0)}_{m+s_0}\geq \beta^{(1)}_i$. Then we put $\theta (\beta^{(0)}_{m+s_0})=\beta^{(1)}_i$.
 \item Assume that we have defined $\theta(\beta^{(0)}_j)$ for $j=p+1,p+2,\cdots,m+s_0$.
 Let  $\beta^{(1)}_k$ be the maximal element of $\beta^{(1)}\setminus\left\{\theta(\beta^{(0)}_{m+s_0}\cup \cdots
 \cup \beta^{(0)}_{p+2} \cup \beta^{(0)}_{p+1}) \right\}$
  such that $\beta^{(0)}_{p}\geq \beta^{(1)}_k$. Then we put $\theta
  (\beta^{(0)}_{p})=\beta^{(1)}_k$.
\end{itemize}
Observe that the standardness of $S_{\bf s}({\ulambda})$ implies that $\theta$
is well-defined. The $2$-tuples $(j,\theta(j))$ such that $\theta (j)\neq j$
are called the {\it pairs} of $S_{\bf s}({\ulambda})$.

\begin{exmp}
Let $e=4$, ${\bf s}=(0,2)$. Then by Prop. \ref{FL}, the bipartition
$\ulambda:=((2,2,1),(3,2))$ is in $\Phi_{4,10}^{(0,2)}$. The ${\bf s}$-symbol
of this bipartition is the following one (where we put $m=4$):
 $$\left(\begin{array}{ccccccc}
  0 &  1 & 2  & 3 & 6 & 8  \\
 0 & 2 & 4 & 5  &
\end{array}
\right)$$ We have $\theta (0)=0$, $\theta (2)=2$, $\theta (4)=3$, $\theta
(5)=1$.

\end{exmp}

\begin{defn} Let $e$ be a positive integer such that $e>1$ and let
$\boldsymbol{s}:=(s_0,s_1)\in{\mathbb{N}^2}$ be such that $s_0\leq s_1$. We
define a map :
$$\Upsilon_{(s_0,s_1)}: \Phi_{e,n}^{(s_0,s_1)} \to \Pi_{2,n}$$
as follows. Let $\ulambda\in\Phi_{e,n}^{(s_0,s_1)}$, let $S_{\bf
s}({\ulambda})$ be the associated ${\bf s}$-symbol. Let $S_{\bf s}'$ be the
symbol obtained from $S_{\bf s}({\ulambda})$ by permuting the pairs in  $S_{\bf
s}({\ulambda})$ and reordering the rows. Let $\umu$ be the bipartition such
that  $S_{\bf s}'=S_{\bf s}({\umu})$. Observe that  $\umu\in\Pi_{2,n}$. Then we
put:
$$\Upsilon_{(s_0,s_1)}(\ulambda)=\umu$$
\end{defn}
\begin{exmp}
Keeping the above example, the symbol $S_{\bf s}'$ is given by
 $$\left(\begin{array}{cccccc}
   0 & 2 & 4  & 5 & 6 &  8  \\
 0 & 1 &    2 &   3 &   &
\end{array}
\right)$$ This is the ${\bf s}$-symbol of  the bipartition
$(\emptyset,(3,2,2,2,1))$.

\end{exmp}

\begin{rem}
 Note that the reversed map  $\Upsilon^{-1}_{(s_0,s_1)}$ can be easily obtained
 as follows.  Let
${\umu}:=(\mu^{(0)},\mu^{(1)})\in
\Upsilon_{(s_0,s_1)}(\Phi_{e,n}^{(s_0,s_1)})$ and let $S_{\bf
s}({\umu})=\left(\begin{array}{cc}
\beta^{(1)}\\
\beta^{(0)}
\end{array}
\right) $ be its ${\boldsymbol
 s}$-symbol. We define an injective map $\tau:\beta^{(0)}\to
 \beta^{(1)}$ such that $\tau(\beta^{(0)}_j)\geq \beta^{(1)}_j$ for all $j\in{\{1,\cdots,m+s_1\}}$ as follows.
 \begin{itemize}
 \item Let $\beta^{(1)}_i$ be the minimal element of $\beta^{(1)}$ such that
 $\beta^{(0)}_{1}\leq \beta^{(1)}_i$. Then we put $\tau (\beta^{(0)}_{1})=\beta^{(1)}_i$.
 \item Assume that we have defined $\theta(\beta^{(0)}_j)$ for $j=1,2,\cdots,p-1$.
 Let  $\beta^{(1)}_k$ be the minimal element of $\beta^{(1)}\setminus\left\{\tau(\beta^{(0)}_{1}\cup \beta^{(0)}_{2} \cup \cdots
 \cup \beta^{(0)}_{p-1}) \right\}$
  such that $\beta^{(0)}_{p}\leq \beta^{(1)}_k$. Then we put $\tau
  (\beta^{(0)}_{p})=\beta^{(1)}_k$.
\end{itemize}
Let $\ulambda$ be the bipartition associated to the ${\bf s}$-symbol obtained
from $S_{\bf s}({\umu})$ by permuting the pairs $(j,\tau(j))$ with $j\neq \tau
(j)$ and reordering the rows. Then we have $\umu=\Upsilon^{-1}_{(s_0,s_1)}
(\ulambda)$.

\end{rem}

\subsection{Main result} We can now state the main theorem of this paper which
gives the explicit description of the bijection
$\Psi_{(s_0,s_1)}^{(s_0,s_1+e)}$.

\begin{thm}\label{main}
Let $e$ be a positive integer such that $e>1$ and let
$\boldsymbol{s}:=(s_0,s_1)\in{\mathbb{N}^2}$ be such that $s_0\leq s_1$ then:
$$\Psi_{(s_0,s_1)}^{(s_0,s_1+e)}=\Upsilon_{(s_0,s_1)}$$
\end{thm}
To prove this theorem, we will need combinatorial properties of the map
$\Upsilon_{(s_0,s_1)}$. Recall that $m\in{\mathbb{N}}$ is such that $m>
\textrm{Max}(r_0-s_0,r_1-s_1)$. For a bipartition $\unu \in\Pi_{2,n}$, let
$S_{\bf s}({\unu})=\left(\begin{array}{cc}
\beta^{(1)}\\
\beta^{(0)}
\end{array}
\right) $ be its ${\bf
 s}$-symbol. Observe that each node $\gamma$ on the {\it border} of $\unu$ (that is
 at the right ends of the Young diagram of $\unu$)
 corresponds to  an element of   $S_{\bf
s} (\unu)$. Indeed, to each node $(a,\nu^{(c)}_a,c)$, we can associate the
element $\beta^{(c)}_a=\nu^{(c)}_a-a+s_c+m$.
 Observe also that:
\begin{itemize}
\item If the number $\beta^{(c)}_a -1$ doesn't occur in  $\beta^{(c)}$ then  $\gamma$ is a removable node of
$\nu$.
\item If the number $\beta^{(c)}_a +1$ doesn't occur in  $\beta^{(c)}$ then we
have an addable node  $\gamma':=(a,\nu_a^{(c)}+1,c)$ in $\nu$.
\item The residue of the node $\gamma$ associated to $\beta^{(c)}_a$ is $\beta^{(c)}_a-m (\textrm{mod
}e)$.

\end{itemize}
In addition, recall that if $\eta=(a,b,c)$ and $\eta'=(a',b',c')$ are two
$i$-nodes of a bipartition, We have $\eta<_{(s_0,s_1)}\eta'$ if and only if :
$$b-a+s_c<b'-a'+s_{c'}\ \textrm{or } \textrm{if}\ b-a+s_c=b'-a'+s_{c'}\textrm{ and }{c}>{c'}.$$
On the other hand, assume that $\eta=(a,b,c)$ and $\eta'=(a',b',c')$ are two
$i$-nodes such that  $\eta<_{(s_0 ,s_1 +e)}\eta'$.
\begin{itemize}
\item if $c=c'=0$ then we have $b-a+s_0<b'-a'+s_{0}$,
\item if $c=c'=1$ then we have $b-a+s_1<b'-a'+s_{1}$,
\item if $c=0$ and $c=1$ then we have $b-a+s_0<b'-a'+s_{1}+e$. Thus we have
$b-a+s_0<b'-a'+s_{1}$ or $b-a+s_0=b'-a'+s_{1}$,
\item if $c=1$ and $c'=0$, we have $b-a+s_1+e<b'-a'+s_{0}$. Thus we have
$b-a+s_1<b'-a'+s_{0}$.
\end{itemize}
Hence, if $\eta=(a,b,c)$ and $\eta'=(a',b',c')$ are two $i$-nodes of a
bipartition, we have $\eta<_{(s_0 ,s_1 +e)}\eta'$ if and only if :
$$b-a+s_c<b'-a'+s_{c'}\ \textrm{or } \textrm{if}\ b-a+s_c=b'-a'+s_{c} \textrm{ and
}{c}<{c'}.$$

\subsection{Proof of Theorem \ref{main}}
This is proved by induction on $n$. If $n=0$ then the result is trivial as
$$\Psi_{(s_0,s_1)}^{(s_0,s_1+e)}({\boldsymbol \emptyset})=\Upsilon_{(s_0,s_1)}=({\boldsymbol\emptyset})={\boldsymbol\emptyset}.$$
Let $n>0$, let $\ulambda:=(\lambda^{(0)},\lambda^{(1)})\in\Phi_{e,n}^{\bf s} $
and let $S_{\bf s}({\ulambda})=\left(\begin{array}{cc}
\beta^{(1)}\\
\beta^{(0)}
\end{array}
\right) $ be its ${\boldsymbol
 s}$-symbol. Let
$\umu=(\mu^{(0)},\mu^{(1)}):=\Upsilon_{(s_0,s_1)}(\lambda^{(0)},\lambda^{(1)})$
and let   $S_{\bf s}({\umu})=\left(\begin{array}{cc}
\alpha^{(1)}\\
\alpha^{(0)}
\end{array}
\right) $ be its ${\bf s}$-symbol.

As in \S \ref{crystal},  we write the sequence of removable and addable
$i$-nodes of $\ulambda$  in increasing order with respect to $<_{(s_0,s_1)}$ :
 $$A_1 A_2 R_3 R_4 A_5 R_6 \cdots A_s$$
 where we write $R_j$ for a removable $i$-node and $A_j$ for an addable
 $i$-node. We delete the occurrences $R_j A_{j+1}$ in this sequence. Then, we obtain a
sequence $\mathfrak{S}$  of removable $i$-nodes  and addable $i$-nodes:
$$A_{j_1} \cdots A_{j_s}R_{i_1} R_{i_2} \cdots$$
 where $j_1<j_2<\cdots<i_1<i_2<\cdots$. The $R_{i_k}$ correspond to the normal $i$-nodes of $\ulambda$ and the leftmost
 one, $R_{i_1}$, is a good $i$-node for $\ulambda$.

Let $R_{i_l}$ be an element of  $\mathfrak{S}$. As explained above, $R_{i_l}$
corresponds to an element $\beta^{(c)}_a$ in  $S_{\bf s}({\ulambda})$. As
$R_{i_l}$ is removable, we have $\beta_{a-1}^{(c)}< \beta_a^{(c)}-1$. We will
associate to this node a removable $i$-node $R_{i_l}'$ in $\umu$. To do this,
we will distinguish several cases. In each case, we give an example of the
symbols $S_{\bf s}({\umu})$ and $S_{\bf s}({\ulambda})$ in which the elements
corresponding to $R_{i_l}$ and $R_{i_l}'   $ are written in bold.

\begin{enumerate}
\item Assume that $c=0$ and that  we have $\theta (\beta_a^{(0)})=\beta_b^{(1)}< \beta_a^{(0)}$ for
 $b\in{\{1,\cdots,m+s_1\}}$. Then to obtain $S_{\bf
s}({\ulambda})$, we have to permute   $\beta_a^{(0)}$ and  $\beta_b^{(1)}$. As
$\beta_b^{(1)}< \beta_a^{(0)}$, the node $R_{i_l}'$ associated to
$\beta^{(0)}_a$ in $\alpha^{(1)}$ is a removable $i$-node (because
$\beta^{(0)}_a -1$ cannot occur in $\alpha^{(1)}$). Note that if we have
$\beta_b^{(1)}=\beta^{(0)}_a-1$, then we have an addable $i$-node $A$ on the
part of $\ulambda$ associated to $\beta^{(1)}_b$ in $\beta^{(1)}$ such that
$A<_{(s_0,s_1)}R_{i_l}  $. In this case, we have an addable $i$-node $A'$ on
the part of $\umu$ associated to $\beta^{(1)}_b$ in $\alpha^{(0)}$ such that
$A'<_{(s_0,s_1+e)}R_{i_l}'  $. \begin{exmp} In the following example, we put
$\beta^{(0)}_a=j$, $\beta^{(0)}_b=j-1$ $\beta_{a-1}^{(0)}=j-2$ and
$\beta_{b-1}^{(1)}=j-3$.
 $$S_{\bf s}({\ulambda})=\left(\begin{array}{ccccc}
    \cdots &  j-3  & j-1 & \cdots & \cdots   \\
  \cdots &    j-2 & {\bf j} & \cdots  &
\end{array}
\right)$$ Then
 $$S_{\bf s}({\umu})=\left(\begin{array}{ccccc}
    \cdots &  j-2  & {\bf j} & \cdots & \cdots   \\
  \cdots &    j-3 &  j-1 & \cdots  &
\end{array}
\right)$$
\end{exmp}
\item Assume that $c=0$ and that  we have $\theta (\beta_a^{(0)})=\beta_b^{(1)}= \beta_a^{(0)}$ for
 $b\in{\{1,\cdots,m+s_1\}}$ and that $\beta_{b-1}^{(1)}<\beta^{(1)}_{b}-1$. In this case, we have a removable
 $i$-node $R$ associated $\beta_b^{(1)}$ in $\beta^{(1)}$. Observe that $R<_{(s_0,s_1)} R_{i_l} $.
  Then to obtain $S_{\bf
s}({\ulambda})$, $\beta_a^{(0)}$ isn't permuted with any elements of
$\beta^{(1)}$. The node $R_{i_l}'$ associated to $\beta^{(1)}_b$ in
$\alpha^{(1)}$ is a removable $i$-node. Note that the removable $i$-node $R'$
associated to $\beta^{(0)}_a$ in $\alpha^{(0)}$ is such that
$R'<_{(s_0,s_1+e)}R_{i_l}'$.
 \begin{exmp} In the following example, we put
$\beta^{(0)}_a=j=\beta^{(0)}_b$ $\beta_{a-1}^{(0)}=j-2$ and
$\beta_{b-1}^{(1)}=j-3$.
 $$S_{\bf s}({\ulambda})=\left(\begin{array}{ccccc}
    \cdots &  j-3  & j & \cdots & \cdots   \\
  \cdots &    j-2 & {\bf j} & \cdots  &
\end{array}
\right)$$ Then
 $$S_{\bf s}({\umu})=\left(\begin{array}{ccccc}
    \cdots &  j-2  & {\bf j} & \cdots & \cdots   \\
  \cdots &    j-3 &  j & \cdots  &
\end{array}
\right)$$
\end{exmp}
\item Assume that $c=0$ and that  we have $\theta (\beta_a^{(0)})=\beta_b^{(1)}= \beta_a^{(0)}$ for
 $b\in{\{1,\cdots,m+s_1\}}$ and that $\beta_{b-1}^{(1)}=\beta^{(1)}_{b}-1$.
  Then to obtain $S_{\bf
s}({\ulambda})$, $\beta_a^{(0)}$ isn't permuted with any elements of
$\beta^{(1)}$. The node $R_{i_l}'$ associated to $\beta^{(0)}_a$ in
$\alpha^{(0)}$ is a removable $i$-node.
\begin{exmp} In the following example, we put $\beta^{(0)}_a=j=\beta^{(0)}_b$,
 $\beta_{a-1}^{(0)}=j-2$ and $\beta_{b-1}^{(1)}=j-1$.
 $$S_{\bf s}({\ulambda})=\left(\begin{array}{ccccc}
    \cdots &  j-1  & j & \cdots & \cdots   \\
  \cdots &    j-2 & {\bf j} & \cdots  &
\end{array}
\right)$$ Then
 $$S_{\bf s}({\umu})=\left(\begin{array}{ccccc}
    \cdots &  j-1  &  j & \cdots & \cdots   \\
  \cdots &    j-2 &  {\bf j} & \cdots  &
\end{array}
\right)$$\end{exmp}
\item Assume that $c=1$ and that  we have $\theta (\beta_b^{(0)})=\beta_a^{(1)}< \beta_b^{(0)}$ for
 a $b\in{\{1,\cdots,m+s_1\}}$.
  Then to obtain $S_{\bf
s}({\ulambda})$, $\beta_b^{(0)}$ must be permuted with $\beta^{(1)}_a$. The
node $R_{i_l}'$ associated to $\beta^{(1)}_a$ in $\alpha^{(0)}$ is a removable
$i$-node.  Note that if we have $\beta_{b-1}^{(0)}=\beta^{(1)}_a-1$, then we
have an addable $i$-node $A$ on the part of $\ulambda$ associated to
$\beta^{(0)}_{b-1}$ such that $ A >_{(s_0,s_1)} R_{i_l}$. In this case, we have
an addable $i$-node $A'$ on the part of $\umu$ associated to
$\beta^{(0)}_{b-1}$ in $\alpha^{(1)}$ such that $A'>_{(s_0,s_1+e)} R_{i_l}'$.
\begin{exmp}
In the following example, we put $\beta^{(1)}_a=j$, $\beta^{(0)}_b=j+1$,
$\beta_{a-1}^{(1)}=j-3$ and $\beta_{b-1}^{(1)}=j-1$.
 $$S_{\bf s}({\ulambda})=\left(\begin{array}{ccccc}
    \cdots &  j-3  & {\bf j} & \cdots & \cdots   \\
  \cdots &    j-1 &  j+1 & \cdots  &
\end{array}
\right)$$ Then
 $$S_{\bf s}({\umu})=\left(\begin{array}{ccccc}
    \cdots &  j-1  &  j+1 & \cdots & \cdots   \\
  \cdots &    j-3 &  {\bf j} & \cdots  &
\end{array}
\right)$$\end{exmp}
\item Assume that $c=1$ and that  we have $\theta (\beta_b^{(0)})=\beta_a^{(1)}= \beta_b^{(0)}$ for
 a $b\in{\{1,\cdots,m+s_1\}}$.  Then to obtain $S_{\bf
s}({\ulambda})$, $\beta_a^{(1)}$ isn't permuted with any elements of
$\beta^{(0)}$. The node $R_{i_l}'$ associated to $\beta^{(0)}_b$ in
$\alpha^{(0)}$ must be  a removable $i$-node. Note that if
$\beta_{b-1}^{(0)}<\beta^{(0)}_b-1$, then the node $R$ associated to
$\beta^{(0)}_b$ in $\beta^{(0)}$ is a removable $i$-node such that
$R>_{(s_0,s_1)}R_{i_l}$. Then, the node $R'$ associated to $\beta^{(1)}_a$ in
$\alpha^{(1)}$ is a removable $i$-node such that $R'>_{(s_0,s_1+e)}R_{i_l}'$.
\begin{exmp}
In the following example, we put $\beta^{(0)}_a=j=\beta^{(0)}_b$,
$\beta_{a-1}^{(0)}=j-2$ and $\beta_{b-1}^{(1)}=j-1$.
 $$S_{\bf s}({\ulambda})=\left(\begin{array}{ccccc}
    \cdots &  j-2  & {\bf j} & \cdots & \cdots   \\
  \cdots &    j-1 &  j & \cdots  &
\end{array}
\right)$$ Then
 $$S_{\bf s}({\umu})=\left(\begin{array}{ccccc}
    \cdots &  j-1  &  j & \cdots & \cdots   \\
  \cdots &    j-2 &  {\bf j} & \cdots  &
\end{array}
\right)$$\end{exmp}
\item Assume that $c=1$ and that  we have $\theta (\beta_b^{(0)})\neq \beta_a^{(1)}$ for
 all $b\in{\{1,\cdots,m+s_0\}}$. Then the node $R_{i_l}'$ associated to $\beta^{(1)}_a$ in $\beta^{(1)}$
  is a removable $i$-node except possibly in the following case: there exists $d\in{\{1,\cdots,m+s_0\}}$ such that
   $\beta^{(0)}_d=\beta^{(1)}_a-1$. In this case, we have an addable $i$-node $A$ in $\beta^{(0)}$ such that
   $A>_{(s_0,s_1)} R_{i_l}$ and there is no removable $i$-node
   between $R_{i_l}$ and $A$ in $\ulambda$ contradicting the fact that
   $R_{i_l}$ is a normal $i$-node.
    \begin{exmp}In the following example, we put $\beta^{(1)}_a=j$,
 $$S_{\bf s}({\ulambda})=\left(\begin{array}{ccccc}
    \cdots &  j-3  & {\bf j} & j+1 & \cdots   \\
  \cdots &    j-2 &  j+3 & \cdots  &
\end{array}
\right)$$ Then
 $$S_{\bf s}({\umu})=\left(\begin{array}{ccccc}
    \cdots &  j-2  &   {\bf j} & j+3 & \cdots   \\
  \cdots &    j-3 &   j+1 & \cdots  &
\end{array}
\right)$$ \end{exmp}

\end{enumerate}

Thus we have associated to each normal $i$-node $R_{i_l}$ in $\ulambda$  a
removable $i$-node $R_{i_l}'$ in $\umu$.

Similarly, one can easily check that if $A_{i_t}$ is an addable $i$-node of
$\ulambda$ in $\mathfrak{S}$, then  we can associate an addable $i$-node
$A_{i_t}'$ in $\umu$ as above. If otherwise,  one can show that there exists  a
removable $i$-node $R$ such that $R<_{(s_0,s_1)} A_{i_t}$ and there is no
addable or removable $i$-node between these two $i$-nodes for the order
$>_{(s_0,s_1)}$. This contradicts the fact that we have deleted all the
occurrences $R_j A_{j+1}$ in the sequence $\mathfrak{S}$.

Hence, we have associated to the sequence $\mathfrak{S}$, a sequence of
removable and addable $i$-nodes of $\umu$ :
$$A_{j_1}' \cdots A_{j_s}'R_{i_1}' R_{i_2}'\cdots $$
 where $j_1<j_2<\cdots<i_1<i_2<\cdots$. Note that by the above observations, this
 sequence is written in increasing order with respect to the order
 $<_{(s_0,s_1+e)}$.

 Now, by the above observations, it is easy to verify that this sequence
 correspond to the sequence $\mathfrak{S}'$ of the removable and addable $i$-nodes of $\umu$,
 written in increasing order with respect to $<_{(s_0,s_{1}+e)}$ and where the
 occurrences $RA$ have been deleted.

The only problem may appear in the following situation.  We have an $i$-node
corresponding to an element $\beta^{(1)}_a$ which is not
 removable, there exists $d\in{\{1,\cdots,m+s_0\}}$ such that
   $\theta(\beta^{(0)}_d)=\beta^{(1)}_a<\beta^{(0)}_d$ and
   $\beta^{(0)}_{d-1}<\beta^{(1)}_a-1$. In this situation, to obtain  $S_{\bf s}({\ulambda})$, we must
   permute $\beta^{(0)}_d$ and $\beta^{(1)}_a$. Moreover, $\beta^{(1)}_{a-1}$ isn't
   permuted with any elements of $\beta^{(0)}$. Thus the $i$-node $R'$
   associated to $\beta^{(1)}_a$ in $\alpha^{(0)}$ must be removable for
   $\umu$.
\begin{exmp} In the following example, we put $\beta^{(1)}_a=j$,
$\beta^{(0)}_d=j+1$ and  $\beta_{d-1}^{(0)}=j-2$.
 $$S_{\bf s}({\ulambda})=\left(\begin{array}{cccccc}
  \cdots &   j-3 &  j-1  &  {\bf j} &\cdots  & \cdots   \\
 \cdots &\cdots &    j-2 &  j+1 & \cdots  &
\end{array}
\right)$$ Then
 $$S_{\bf s}({\umu})=\left(\begin{array}{cccccc}
  \cdots  &j-2 &  j-1  &    j+1 & \cdots & \cdots   \\
 \cdots & \cdots &    j-3 &   {\bf j} & \cdots  &
\end{array}
\right)$$
\end{exmp}
 Note that  in this case, we have an addable $i$-node $A '$ on the part of
$\mu^{(1)}$ associated to $\beta^{(1)}_{a-1}$ such that $ A'
>_{(s_0,s_1+e)}R '$. Thus, to obtain $\mathfrak{S}'$, the occurrence $R'A'$
must be deleted.

Now the leftmost removable $i$-node $R_{i_1}$ in $\mathfrak{S}$ is a good
$i$-node for $\ulambda$ (with respect to $<_{(s_0,s_{1})}$) and the above
discussion shows that this corresponds to a removable $i$-node $R_{i_1} '$ in
$\mathfrak{S}'$ which must be a good $i$-node for $\umu$ (with respect to
$<_{(s_0,s_{1}+e)}$).

Finally, let $\unu$ be the bipartition obtained by removing $R_{i_1}$ from
$\ulambda$. Note that in case $(2)$ above, the normal $i$-node $R_{i_l}$ cannot
be a good $i$-node. Indeed, we have a removable $i$-node $R$ such that
$R<_{(s_0,s_1)} R_{i_l}$ and no addable $i$-node between these two nodes. Hence
$R$ is a normal $i$-node such that $R<_{(s_0,s_1)} R_{i_l}$ and thus $R_{i_l}$
is not a good $i$-node.
 Studying the
other cases above, one can verify that $\Upsilon_{(s_0,s_1)}(\unu)$ is the
bipartition obtained by removing $R_{i_1}'$ from $\umu$. This concludes the
proof of the main Theorem.
\subsection{Example} In this subsection, we give an example for the
computation of the bijection $\Psi_{(s_0,s_1)}^{(s_0,s_1+e)}$. We put ${\bf
s}=(0,1)$, $e=4$ and we consider the bipartition
$\ulambda:=((8),(4))\in{\Phi_{4,12}^{(0,1)}}$ (this is a FLOTW bipartition, see
Prop. \ref{FL}). We compute the $(0,1)$-symbol of $((8),(4))$ (with $m=2$).
$$S_{(0,1)}((8),(4) )=\left(\begin{array}{ccc}
    0 & 1  & 6   \\
  0 &    9  &
\end{array}\right)$$
Then, the injection $\theta:\{0,9\}\to \{0,1,6\}$ is such that $\theta (0)=0$
and  $\theta (9)=1$. Thus, the $(0,1)$-symbol of $\Psi_{(0,1)}^{(0,5)}((8),(4))
$ is:
$$\left(\begin{array}{ccc}
    0 & 1  & 9   \\
  0 &    6  &
\end{array}\right)$$
Thus, we  have $\Psi_{(0,1)}^{(0,5)}((8),(4))=((5),(7))$. We now want to find
$\Psi_{(0,5)}^{(0,9)}((5),(7))$. The $(0,9)$-symbol of $((5),(7))$ is :
$$S_{(0,5)}((5),(7)))=\left(\begin{array}{ccccccc}
    0 & 1  & 2& 3& 4 &5 &13  \\
  0 &    6 &
\end{array}\right)$$
Then, the injection $\theta:\{0,6\}\to \{0,1,2,3,4,5,13\}$ is such that $\theta
(0)=0$ and  $\theta (6)=5$. Thus, the $(0,5)$-symbol of
$\Psi_{(0,5)}^{(0,9)}((8),(4)) $ is:
$$\left(\begin{array}{ccccccc}
    0 & 1  & 2&3&4&6&13  \\
  0 &    5  &
\end{array}\right)$$
Thus, we  have $\Psi_{(0,5)}^{(0,9)}((8),(4))=((4),(7,1))$. We want now to find
the Uglov bipartition $\Psi_{(0,9)}^{(0,13)}((4),(7,1))$. The $(0,9)$-symbol of
$((4),(7,1))$ is :
$$S_{(0,9)}((4),(7,1)))=\left(\begin{array}{ccccccccccc}
    0 & 1  & 2& 3& 4 &5&6&7&8&9 &17  \\
  0 &    6 &
\end{array}\right)$$
Then, the injection $\theta:\{0,6\}\to \{0,1,2,3,4,5,6,7,8,9,17\}$ is such that
$\theta (0)=0$ and  $\theta (6)=6$.

Thus, the $(0,9)$-symbol of $\Psi_{(0,5)}^{(0,9)}((8),(4)) $ is
$S_{(0,9)}((4),(7,1))$. Hence we  have
$\Psi_{(0,9)}^{(0,13)}((4),(7,1))=((4),(7,1))$. Now, we have
$\Phi_{4,12}^{(0,13)}=\Phi_{4,12}^{(0,1)_{+}}$ because $13-0>n-1$.

\section{Relation with results of Leclerc and Miyachi}

 Following the works of
Leclerc and Miyachi and using the above results, it is possible to describe the
bijection $\Psi_{(s_0,s_1)}^{(s_0,s_1+e)}$ using the theory of canonical basis
for $\mathcal{U}_v(\mathfrak{sl}_{\infty})$-modules. We first recall the
results of \cite{LM}.

Let $\mathcal{U}_v(\mathfrak{sl}_{\infty})$ be the quantum algebra associated
to the doubly infinite diagram of type $A_{\infty}$. The fundamental weights
are denoted by $\Lambda_i$ with $i\in{\mathbb{Z}}$. Let ${\bf
s}:=(s_0,s_1)\in{\mathbb{Z}^2}$ with $s_0\leq s_1$ and let $L_{{\bf s}}$ be the
irreducible highest weight module with highest weight
$\Lambda_{s_0}+\Lambda_{s_1}$. Then the theory of Kashiwara and Lusztig
provides us with a canonical basis for $L_{{\bf s}}$. This basis is naturally
labelled by the vertices of the associated crystal graph which may be
constructed as in \S \ref{crystal}. It is easy to see that the class of
bipartitions which label this graph is given by :
$$\Phi_{\infty,n}^{{\bf s}}=\left\{ \ulambda=(\lambda^{(0)},\lambda^{(1)})\in{\Pi_{2,n}}\
|\  \lambda_i^{(0)}\geq{\lambda^{(1)}_{i+s_{1}-s_{0}}},\ i=1,2,3,\cdots\right\}
$$ with $n\geq 0$. Thus, if $\ulambda\in{\Phi_{\infty,n}^{{\bf s}}}$, the
associated element of the canonical basis is given by :
$$b(\ulambda)=\sum_{\umu\in{\Pi_{2,n}}} c_{\ulambda,\umu}(v) \umu$$
with $c_{\ulambda,\ulambda}(v)=1$ and
$c_{\ulambda,\umu}(v)\in{v\mathbb{Z}[v]}$ if $\umu\neq \ulambda$.

Let $\ulambda:=(\lambda^{(0)},\lambda^{(1)})\in{\Phi_{\infty,n}^{{\bf s}}}$ and
let $S_{{\bf s}}(\ulambda)$ be its associated symbol. By the above
characterization of $\Phi_{\infty,n}^{{\bf s}}$, this symbol is standard. Let
$p$ be the numbers of pairs in this symbol (see \S \ref{symbol}) and let
$C(\ulambda)$ be the set of bipartitions $\umu$ of $n$ such that $S_{{\bf
s}}(\umu)$ is obtained from $S_{{\bf s}}(\ulambda)$ by permuting some pairs in
$S_{{\bf s}}(\ulambda)$ and reordering the rows. For $\umu\in{C(\ulambda)}$, we
denote by $l(\umu)$ the number of pairs permuted in $S_{{\bf s}}(\ulambda)$ to
obtain $S_{{\bf s}}(\umu)$. In particular, we have $l(\ulambda)=0$.  Then, the
following result gives an explicit description of the canonical basis.
\begin{thm}[Leclerc-Miyachi \protect{\cite[Theorem 3]{LM}}]\label{lm}  Let $\ulambda\in{\Phi_{\infty,n}^{{\bf
s}}}$ and let $b(\ulambda)$ be the associated element of the canonical basis of
$L_{{\bf s}}$. Then, we have :
$$b(\ulambda)=\sum_{\umu\in{C(\ulambda)}} v^{l(\umu)} \umu.$$
\end{thm}
Now, let $e$ be a positive integer such that $e>1$ and let ${\bf
s}:=(s_0,s_1)\in{\mathbb{Z}^2}$ with $s_0\leq s_1$. Let
$\ulambda\in{\Phi_{e,n}^{{\bf s}}}$. By Prop. \ref{std}, we have
$\ulambda\in{\Phi_{\infty,n}^{{\bf s}}}$. Thus $\ulambda$ is labelling the
element of the canonical basis of the irreducible highest weight module
$L_{{\bf s}}$ with highest weight $\Lambda_{s_0}+\Lambda_{s_1}$. Hence Theorem
\ref{main} together with Theorem \ref{lm} yields the following remarkable
property:
\begin{thm}  Let $\ulambda\in{\Phi_{e,n}^{{\bf s}}}$. Then
we have $\Psi_{(s_0,s_1)}^{(s_0,s_1+e)}(\ulambda)=\umu$ if and only if  the
degree of $c_{\ulambda,\umu}(v)$ is maximal in $b(\ulambda)$.
\end{thm}
It could be interesting to obtain a non combinatorial proof of the above
theorem which shows why the bijections $\Psi_{(s_0 ,s_1 )}^{(s_0 ',s_1 ')}$ is
controlled by the canonical basis of irreducible
$\mathcal{U}_v(\mathfrak{sl}_{\infty})$-modules.

Another open problem would be to obtain similar statements for the irreducible
highest weight  $\mathcal{U}_v(\mathfrak{sl}_{\infty})$-modules  of level
$l>2$. In this case, relations between the sets of Uglov multipartitions and
the representation theory of Ariki-Koike algebras have been established in
\cite{J2}.


\end{document}